\documentclass[12pt]{article}
\usepackage[ansinew]{inputenc}
\usepackage[pdftex]{graphicx}
\usepackage{amssymb,amsmath,amscd}
\usepackage{dsfont}
\usepackage{verbatim}
\usepackage{amsmath} 
\usepackage{graphicx}
\usepackage{pstricks}
\usepackage{enumitem}
\usepackage{ntheorem}
\usepackage{amssymb}
\usepackage{float,caption}

\newtheorem{de}{Definition}[section] 
\newtheorem{nott}{Notation}

\newtheorem{theo}{Theorem}
\newtheorem{prop}[theo]{Proposition}  
\newtheorem{lem}{Lemma}[section]
 
\newtheorem{conj}{Conjecture}[section]

\setenumerate[1]{font=\bfseries,label=\Roman*.}
\setenumerate[2]{font=\itshape,label=\arabic*)}

\author{Yashar Memarian
}       
\title{On a Correlation Inequality for Cauchy Type Measures\\
       }
\date{\today}

\begin{document}
\maketitle

\begin{abstract}
In this paper we present a correlation inequality with respect to Cauchy type measures. To prove our inequality, we transport the problem onto the Riemannian sphere then state and solve some special cases for a spherical correlation problem. This method, as we shall explain, opens up a new class of interesting problems related to correlation type inequalities.
\end{abstract}

\section{Introduction}

Perhaps the most important correlation-related problem is the \emph{famous} Gaussian Correlation Conjecture. The standard Gaussian measure (denoted by $\gamma_n$) of any measurable subset $A\subseteq \mathbb{R}^n$ is defined by
\begin{eqnarray*}
\gamma_n(A)=\frac{1}{(2\pi)^{n/2}}\int_{A}e^{-\vert x\vert^2/2} dx.
\end{eqnarray*}
A general mean zero Gaussian measure, $\mu_n$, defined on $\mathbb{R}^n$ is a linear image of the standard Gaussian measure. The Gaussian Correlation Conjecture is formulated as follows :
\begin{conj} \label{1}
For any $n\geq 1$, if $\mu$ is a mean zero, Gaussian measure on $\mathbb{R}^n$, then for $K, M$, convex closed subsets of $\mathbb{R}^n$ which are symmetric about the origin, we have
\begin{eqnarray*}
\mu_n(K\cap M)\geq \mu_n(K)\mu_n(M).
\end{eqnarray*}
\end{conj}
For some background on the above conjecture, a less general form of the Gaussian Correlation Conjecture first appeared in $1955$ in \cite{dunnet}. The general setting appeared a few years later in $1972$ by S. Das Gupta, M.L. Eaton, I. Olkin, M. Perlman, L.J. Savage and M. Sobel in \cite{dasgu}. I won't go into details about what is known or unknown about Conjecture \ref{1}, in this paper, I am more interested about a correlation inequality with respect to another measure. In the past few years, there has been a bit of research about correlation inequalities with respect to measures other than the Gaussian (see for example \cite{lewi1} and \cite{lewi2}). People began to wonder (since it is hard to prove (or disprove) Conjecture \ref{1})if perhaps dealing with other measures could be easier. For example, in \cite{figa}, the authors prove some sharp correlation type inequalities for rotationally invariant measures (with a decay condition on the density function) in $\mathbb{R}^2$. One should recall that Conjecture \ref{1} is proved in the two dimensional case. This paper concerns the following:
\begin{conj} \label{pri}
For any $n\geq 1$, for every two symmetric convex sets $K,M\subset \mathbb{R}^n$, we have
\begin{eqnarray*}
\nu_n(K\cap M)\geq \nu_n(K)\nu_n(M),
\end{eqnarray*}
where 
\begin{eqnarray*}
\nu_n=C.(1+\vert x\vert^2)^{-\frac{n+1}{2}},
\end{eqnarray*}
for $C$ the normalisation constant.
\end{conj}

In dimension $2$, the results proven in \cite{figa} give a positive answer to Conjecture \ref{pri}. Therefore there is no need to discuss this conjecture for this case. Here we discuss several particular cases in higher dimensions for which Conjecture \ref{pri} holds. We also give some ideas on possible proves for this conjecture. The method used to approach Conjecture \ref{pri} is purely geometric. We first announce a spherical correlation conjecture on the canonical Riemannian sphere. Later on, we show that the spherical correlation conjecture is equivalent to Conjecture \ref{pri} using a projective mapping of the Euclidean space to a ball of radius $\pi/2$ of the sphere. Under the projective mapping, every straight line of the Euclidean space is mapped to a geodesic of the sphere. Therefore any convex set of the Euclidean space is mapped to a convex set of the sphere. Morever, the push-forward of the Cauchy measure $\nu_n$ with density $C.(1+\vert x\vert^2)^{-\frac{n+1}{2}}$ (where $C>0$ is a normalisation constant) with respect to the Lebesgue measure, is mapped to the normalised canonical Riemannian measure of the sphere (or, better said \emph{ball}). The spherical correlation will be the subject of Section $3$ of this paper. In Section $5$, we define the projective mapping connecting Conjecture \ref{pri} to the spherical one. In the final section, we discuss possible ways of proving the spherical correlation conjecture.
\section{Acknowledgement}
I am grateful to Michel Ledoux and Franck Barthe for their useful remarks concerning this project.

\section{A Correlation Conjecture on the Sphere}
In this section, we present a correlation conjecture on the canonical Riemannian sphere. Let $\mathbb{S}^n$ be the canonical Riemannian sphere. Fix a hemi-sphere $\mathbb{S}^n_{+} \subset \mathbb{S}^n$. We denote the center of the hemi-sphere by $o$ and recall that $\mathbb{S}^n_{+}=B(o,\pi/2)$ where $B(o,\pi/2)$ is a spherical ball of radius $\pi/2$ centered at the point $o$. We denote the Riemannian volume of the sphere by $vol_n$.
\begin{de}
An open set $S\subset B(o,\pi/2)\subset \mathbb{S}^n$ is convex if it is geodesically convex with respect to the canonical Riemannian geometry of the sphere.  A convex set $S$ is centrally symmetric around a point $x\in Int(S)$ if for any geodesic segment $\sigma$ containing the point $x$
\begin{eqnarray*}
l([x,x^{+}_b])=l([x,x^{-}_b]).
\end{eqnarray*}
Where $\sigma\cap \partial S=\{x^{+}_b,x^{-}_b\}$ and $l(.)$ stands for the length (which is understood as the length related to the Riemannian structure of the sphere).
\end{de}
For our future purposes it's best to set the following :
\begin{nott}
Let $X$ be a general metric-measure space and $Y\subset X$. Then
\begin{eqnarray*}
Y+\varepsilon =\{x \in X \vert \hspace{0.5mm} d(x,Y)\leq \varepsilon\},
\end{eqnarray*}
where $d(.,.)$ stands for the metric of $X$ and $d(x,Y)=inf_{y\in Y} d(x,y)$.

\end{nott}

We shall need to remind (and define) a fairly known operation on subsets of the sphere which will come useful in the next section :
\begin{de}[Double Suspension] \label{cone}
Let $X_k\subset\mathbb{S}^k_{o}$ be a $k$-dimensional symmetric convex set containing the point $o$. Let $B_{n-k}^{\perp}(o,\pi/2)$ be the $(n-k)$-dimensional ball \emph{orthogonal} to $\mathbb{S}^k_{o}$ containing $X_k$. Let $\mathbb{S}^{n-k-1}=\partial B_{n-k}^{\perp}(o,\pi/2)$. Define $X_n=X_k *\mathbb{S}^{n-k-1}$ to be the $n$-dimensional symmetric convex set which contains $X_k$ and all the geodesics orthogonal to $X_k$ joining $\mathbb{S}^{n-k-1}$.
\end{de}

\emph{Remark}:

This operation for $k=1$ defines a double cone over a $(n-1)$-dimensional symmetric convex set $X_{n-1}$ and hence generalises it to higher dimensions.

We are now ready to announce the spherical correlation conjecture :

\begin{conj} \label{princ}
Let $K_1$ and $K_2$ be two geodesically convex (spherical) bodies contained in the hemi-sphere $B(o,\pi/2)\subset \mathbb{S}^n$. Additionaly, $K_1$ and $K_2$ are both centrally symmetric around the point $o$. Then
\begin{eqnarray*}
vol_n(K_1\cap K_2).vol_n(B(o,\pi/2))\geq vol_n(K_1).vol_n(K_2).
\end{eqnarray*}
\end{conj}

Conjecture \ref{princ}(in its most general form) is open. However, in the next section we examine a few important specific cases for which there is a positive answer to this conjecture.

\section{A Few Special Cases in Conjecture \ref{princ}}

The aim of this section is to prove the following :
\begin{theo} \label{dmain}
For the following different cases, the result of conjecture \ref{princ} holds :
\begin{itemize}
\item One set is a spherical ball, i.e. $K_i=B(o,r)$ for any $r>0$.
\item Let $0<\varepsilon\leq \pi/2$ and a $1\leq k\leq (n-1)$. Let $Y$ be a spherical tube of width $\varepsilon$, i.e. $Y=\mathbb{S}^k+\varepsilon$ and let $X_n=X_k*\mathbb{S}^{n-k-1}$ as in definition \ref{cone}. Then $X_n$ and $Y$ satisfy the equality case in conjecture \ref{princ}.
\item If an integer $N$ exists such that for every $n\geq N$ conjecture \ref{princ} is true, then the conjecture is true for every dimension.
\item If $X$ is an arbitrary symmetric spherical convex set and $Y$ is any tube, i.e. $Y=\mathbb{S}^k+\varepsilon$, then conjecture \ref{princ} holds for $X$ and $Y$.
\end{itemize}
\end{theo}

Next sections concern the proof of Theorem \ref{dmain}.

\subsection{One Set Is a Spherical Ball}

Let $K_1=B(o,r)$ for $0<r\leq \pi/2$ and $K_2$ an arbitrary convex set (containing $o$).

In this case, the proof of Conjecture \ref{princ} follows immediately applying the following version of the Bishop-Gromov Inequality :

\begin{lem} \label{bg}
For all open convex sets $S\subset B(o,\pi/2)$ and all $x\in S$, the function
\begin{eqnarray*}
\frac{vol_n(S\cap B(x,r))}{vol_n(B(x,r))}
\end{eqnarray*}
is a non-increasing function of $r$.
\end{lem}
\begin{flushright}
$\Box$
\end{flushright}

\emph{Remark}:
\begin{itemize}
\item Note that $S$ is any convex set and not necessarely symmetric.
\item Lemma \ref{bg} is not sharp. There are some non-convex sets in existance for which this Lemma still holds. Indeed for $S$ being a waist of length $r$ around a hyper-sphere i.e. $\mathbb{S}^{n-1}+r$ and for all the tubes $\mathbb{S}^{k}+r$ ($1\leq k\leq n-2$), Lemma \ref{bg} still holds.
\end{itemize}

\subsection{The Equality Case in Conjecture \ref{princ}}

The second natural question to ask is whether it is possible to classify the sets (or at least some class of sets) for which the equality holds in Conjecture \ref{princ}. Normally in a correlation type problem, studying the equality cases is as hard as solving the original general inequality. But, as we shall see in this section, for the spherical correlation inequality, one has a nice characterisation of the equality cases. Of course, one obvious equality case is when one set is the whole half-sphere i.e. $B(o,\pi/2)$. We shall examine a less obvious class of examples for which the equality in Conjecture \ref{princ} holds. Let $X_k$ and $X_n$ be defined as in definition \ref{cone}. One can verify in a straightforward way the following :
\begin{lem} \label{formul}
Let $X_k$ and $X_n=X_n=X_k *\mathbb{S}^{n-k-1}$ be as defined above. Then
\begin{eqnarray*}
\frac{vol_n(X_n)}{vol_n(\mathbb{S}^n_{+})}=\frac{vol_k(X_k)}{vol_k(\mathbb{S}^k_{+})}.
\end{eqnarray*}
\end{lem}
Let $\varepsilon>0$ and let $Y=\mathbb{S}^k+\varepsilon$ be a tube of radius $\varepsilon$. Choose $\mathbb{S}^k$ such that it contains $X_k$. Then
\begin{prop} \label{equal}
For $X_n$ and $Y$ defined as above, we have
\begin{eqnarray*}
vol_n(\mathbb{S}^n_{+}).vol_n(X_n\cap Y)=vol_n(X_n).vol_n(Y).
\end{eqnarray*}
\end{prop}
We recall that the set $Y$ is not a convex set, and it is unclear whether two symmetric convex sets, both different from $B(o,\pi/2)$ exist, such that they satisfy the equality case of Conjecture \ref{princ}.

\emph{Proof of Proposition \ref{equal}}:

Remark that 
\begin{eqnarray*}
vol_n(X_n\cap Y)=(\int_{0}^{\varepsilon}\cos(t)^{k}\sin(t)^{n-k-1}dt).vol_k(X_k).
\end{eqnarray*}
Hence
\begin{eqnarray*}
vol_n(\mathbb{S}^n_{+}).vol_n(X_n\cap Y)&=& vol_n(\mathbb{S}^n_{+}).(\int_{0}^{\varepsilon}\cos(t)^{k}\sin(t)^{n-k-1}dt).vol_k(X_k)\\
                                        &=& (\int_{0}^{\varepsilon}\cos(t)^{k}\sin(t)^{n-k-1}dt). vol_k(\mathbb{S}^k_{+}).vol_n(X_n)\\
                                        &=& vol_n(Y).vol_n(X_n).
\end{eqnarray*}
This ends the proof of Proposition \ref{equal}.

\begin{flushright}
$\Box$
\end{flushright}

\subsection{High Dimensions Imply All Dimensions}

Here we show a simple yet useful result regarding Conjecture \ref{princ}. Roughly speaking, if one can prove the spherical correlation Conjecture for high dimensional spheres, then the conjecture holds in all dimensions. More precisely:
\begin{lem} \label{high}
Suppose an integer $N\in \mathbb{N}$ exists such that for every $n\geq N$-dimensional sphere, Conjecture \ref{princ} holds. Then Conjecture \ref{princ} holds for all $n$.
\end{lem}

\emph{Proof of Lemma \ref{high}}:

Let $k\leq N$ and let $K_1, K_2\subset \mathbb{S}^k_{+}$ two symmetric convex sets around $o \in \mathbb{S}^k_{+}$. See $\mathbb{S}^k_{+}\subset \mathbb{S}^{N}_{+}$ as a $k$-dimensional totally geodesic sub-sphere. Let $X_i=K_i*\mathbb{S}^{N-k-1}$ for $i=1,2$ and $Y=(K_1\cap K_2)*\mathbb{S}^{N-k-1}$. Remark that $X_1, X_2, Y$ are symmetric convex sets in $\mathbb{S}^N$. By the assumptions of Lemma \ref{high} and by the result of Lemma \ref{formul}, we have 
\begin{eqnarray*}
\frac{vol_k(K_1\cap K_2)}{vol_k(\mathbb{S}^{k}_{+})} &=& \frac{vol_N(Y)}{vol_N(\mathbb{S}^{N}_{+})}\\
                                                     &\geq& \frac{vol_N(X_1\cap X_2)}{vol_N(\mathbb{S}^{N}_{+})}\\
                                                     &\geq& \frac{vol_N(X_1)}{vol_N(\mathbb{S}^{N}_{+})}.\frac{vol_N(X_2)}{vol_N(\mathbb{S}^{N}_{+})}\\
                                                     &=& \frac{vol_k(K_1)}{vol_k(\mathbb{S}^{k}_{+})}.\frac{vol_k(K_2)}{vol_k(\mathbb{S}^{k}_{+})}.
\end{eqnarray*}
This ends the proof of Lemma \ref{high}.

\begin{flushright}
$\Box$
\end{flushright}

Lemma \ref{high} is very useful since it is now adequate to prove the spherical correlation conjecture for sufficiently high dimensional spheres.




\subsection{General Symmetric Convex Body and Tubes}

Here we shall examine a harder case for Conjecture \ref{princ}. We assume that none of the convex sets contain the other. The proof of Conjecture \ref{princ} for the opposite case is trivial.

We first need a bit of background on the metric invariant \emph{waist} and an important class of measures on the sphere called $\sin^k$-concave measures.

The waist of a general mm-space is defined in \cite{grwst}. Let $X$ be a mm (metric-measure)-space of dimension $n$. For every $1\leq k\leq n$ the $k$-waist of $X$, denoted by $wst_k(X)$ is the infimum of numbers $r\geq 0$ such that for every family of $k$-cycles (or relative cycles) parametrised by a $n-k$ dimensional $\mathbb{Z}_2$-topological manifold and genereting the fundamental $\mathbb{Z}_2$-homology class of the space of $k$-cycles, the $k$-volume of every cycle is at most equal to $r$. The waist of the canonical sphere is sharply estimated in \cite{grwst} and \cite{memwst}. 

A convex subset $X$ of the canonical Riemannian sphere has sectional curvature everywhere (on its regular part) at least equal to $1$, i.e. $sec(X)\geq 1$. In his recent paper \cite{groexp}, Gromov, by reviewing deeply the ideas of F.Almgren in \cite{almg} and the Heintze-Karcher Volume Comparison Theorem, gives a sharp estimate of the waist of Riemannian manifolds with sectional curvature at least equal to $1$. (The number one can also be relaxed to a $\kappa>0$). More precisely:

\begin{theo} \label{wt}
Let $X$ be a compact connected Riemannian manifold (with a possibly non-empty quasi-regular convex boundary) such that $sec(X)\geq 1$. Let $f:X\to \mathbb{R}^k$ be a smooth map. Then there exists a $z\in \mathbb{R}^k$ such that:
\begin{eqnarray*}
\frac{vol_{n-k}(f^{-1}(z))}{vol_n(X)}\geq \frac{vol_{n-k}(\mathbb{S}^{n-k})}{vol_n(\mathbb{S}^n)}.
\end{eqnarray*}
Where $vol_{n-k}$ stands for the Riemannian volume (or equivalently the Hausdorff measure) in dimension $(n-k)$. 
\end{theo}

For a proof of this theorem one can see \cite{groexp}, or \cite{mempos} where the present author provides a detailed proof of Theorem \ref{wt}.

We recall that a compact $n$-dimensional rectifiable set $X_1\subset X_2$ is called \emph{quasi-regular} if the complementary of the subset of regular parts of $X_1$ has measure zero in $X_1$, and for almost all the points of $X_2$, the distance function $d_x(y)=d(x,y)$ (where the distance is with respect to $X_2$) has its minimum (in $X_1$) at a regular point of $X_1$.

Here, since we are dealing with convex subsets of the sphere, we want to apply Theorem \ref{wt}. Our convex sets may very well have some singularities on the boundary, but since it is enough to use spherical polytopes (which are the interior of the intersection of a finite number of hyperspheres) then we fall under the assumptions of Theorem \ref{wt}. This gives us the desired lower bound for the waist of the convex sets $K_1$, $K_2$ and $K_1\cap K_2$.  

Take a parallel family of hyperspheres $\{\mathbb{S}^{n-1}_t\}_{t\in I}$ which sweep out the convex set $K_1$ (the notion of \emph{parallel} is well defined on symmetric Riemannian manifolds as "one can define two hypersurfaces as parallel if they have parallel second fundamental forms"). $I$ is an interval of $\mathbb{R}$. Since $K_1$ is centrally symmetric, then for every $t\in I$ we have
\begin{eqnarray*}
vol_{n-1}(\mathbb{S}^{n-1}_{o}\cap K_1)\geq vol_{n-1}(\mathbb{S}^{n-1}_t\cap K_1),
\end{eqnarray*}
where $\mathbb{S}^{n-1}_{o}$ is the only hypersphere in this family that contains the point $o$.

Indeed if this is not the case, then there is a $t_0\in I$ such that for every $t\in I$ 
\begin{eqnarray*}
vol_{n-1}(\mathbb{S}^{n-1}_{t_0}\cap K_1)\geq vol_{n-1}(\mathbb{S}^{n-1}_t\cap K_1),
\end{eqnarray*}
and $\mathbb{S}^{n-1}_{t_0}$ does not contain the point $o$. Then by the symmetry of $K_1$ there is another $t'_0\neq t_0$ such that
\begin{eqnarray*}
vol_{n-1}(\mathbb{S}^{n-1}_{t_0}\cap K_1)=vol_{n-1}(\mathbb{S}^{n-1}_{t'_0}\cap K_1)
\end{eqnarray*}
and by convexity of $K_1$ this means that there is a hypersphere $\mathbb{S}^{n-1}_t$ between $\mathbb{S}^{n-1}_{t_{0}}$ and $\mathbb{S}^{n-1}_{t'_0}$ such that
\begin{eqnarray*}
vol_{n-1}(\mathbb{S}^{n-1}_{t}\cap K_1)\geq vol_{n-1}(\mathbb{S}^{n-1}_{t_0}\cap K_1)=vol_{n-1}(\mathbb{S}^{n-1}_{t'_0}\cap K_1),
\end{eqnarray*}
and this is a contradiction. This argument, combined with Theorem \ref{wt}, shows that for any hypersphere $\mathbb{S}^{n-1}_{o}$ which contains the point $o$ we have :
\begin{eqnarray*}
\frac{vol_{n-1}(\mathbb{S}^{n-1}_{o}\cap K_1)}{vol_n(K_1)}\geq \frac{vol_{n-1}(\mathbb{S}^{n-1}_{o}\cap B(o,\pi/2))}{vol_n(B(o,\pi/2))}.
\end{eqnarray*}

Unfortunately, the above waist inequality is not enough for our purpose. We need to extend the waist inequality for the volume of the $\varepsilon$-neighborhood of sections $\mathbb{S}^{n-1}_{o}\cap K_1$. To do this, we require some information and tools about $\sin^k$-measures and functions:

\begin{de}[$\sin$-concave functions]
A real non-negative function $f$ defined on an interval of length less than $2\pi$ is called $\sin$-concave, if, when transported by a unit speed parametrisation of the unit circle, it can be extended to a $1$-homogeneous and concave function on a convex cone of $\mathbb{R}^2$.
\end{de}

\begin{de}[$\sin^k$-concave functions]
A non-negative real function $f$ is called $\sin^k$-concave if the function $f^{\frac{1}{k}}$ is $\sin$-concave.
\end{de}

One can use the following lemma as a definition for $\sin^k$-concave functions:
\begin{lem}
A real non-negative function defined on an interval of length less than $\pi$ is $\sin^k$-concave if for every $0<\alpha<1$ and for all $x_1,x_2\in I$ we have
\begin{eqnarray*}
f^{1/k}(\alpha x_1+(1-\alpha)x_2)\geq (\frac{\sin(\alpha\vert x_2-x_1\vert)}{\sin(\vert x_2-x_1\vert)})f(x_1)^{1/k}+(\frac{\sin((1-\alpha)\vert x_2-x_1\vert)}{\sin(\vert x_2-x_1\vert)}f(x_2)^{1/k}.
\end{eqnarray*}
Particularly if $\alpha=\frac{1}{2}$ we have
\begin{eqnarray*}
f^{1/k}(\frac{x_1+x_2}{2})\geq \frac{f^{1/k}(x_1)+f^{1/k}(x_2)}{2\cos(\frac{\vert x_2-x_1\vert}{2})}.
\end{eqnarray*}
\end{lem}

The following important lemma is proved in \cite{memwst}:
\begin{lem} \label{grmem}
Let $\mu=f.dvol_{k}$ be a measure with a $\sin^{n-k}$-concave density with respect to the $k$-dimensional Riemannian volume of $\mathbb{S}^k$. Let the measure $\mu$ be supported on a $k$-dimensional convex subset $S\subseteq \mathbb{S}^k$. Let $o\in S$ be the point where $f$ attains its maximum, then :
\begin{eqnarray*}
\frac{\int_{B(o,\varepsilon)\cap S}f(x) dvol(x)}{\int_{S}f(x) dvol(x)} &\geq& \frac{\int_{0}^{\varepsilon}\cos^{n-k}(t)\sin^{k-1}(t)dt}{\int_{0}^{\pi/2}\cos^{n-k}(t)\sin^{k-1}(t) dt} \\
                                                                       &=&\frac{vol_n(\mathbb{S}^{n-k}+\varepsilon)}{vol_n(\mathbb{S}^n)}.
\end{eqnarray*}
\end{lem}

The following lemma is a simplified version of a spherical Brunn Theorem:
\begin{lem} \label{brun}
Let $S\subset \mathbb{S}^n$ be a symmetric convex set with respect to $o \in S$. Let $\mathbb{S}^k_{0}\subset \mathbb{S}^n$ be a $k$-dimensional totally geodesic sub-sphere containing the point $o$. Let $\mathbb{S}^k\subset\mathbb{S}^{k+1}_{o}$ be a hyper-sphere in a $k+1$-dimensional sphere containing $o$. Let $p:S\to \sigma$ be the orthogonal projection of $S$ onto $\sigma$. Then  the push-forward of the Riemannian volume has a $\sin^{n-1}$-concave density with respect of $dt$ the canonical measure of the segment $\sigma$.
\end{lem}

We are ready to prove a special case for which Conjecture \ref{princ} holds:

\begin{prop} \label{import}
For every $\varepsilon>0$ and every $\mathbb{S}^{n-1}_{o}$ we have :
\begin{eqnarray*}
\frac{vol_n((\mathbb{S}^{n-1}_{o}+\varepsilon)\cap K_1)}{vol_n(K_1)}\geq \frac{vol_n(\mathbb{S}^{n-1}_{o}+\varepsilon)}{vol_n(\mathbb{S}^n)}.
\end{eqnarray*}
\end{prop}

\emph{Proof of Proposition \ref{import}}:

Let $\varepsilon>0$ and $\mathbb{S}^{n-1}_{o}$ be fixed. Let $\sigma$ be a geodesic orthogonal to $\mathbb{S}^{n-1}_{o}$. Project orthogonally $K_1$ onto $\sigma$. Applying Lemma \ref{brun}, the density of the push-forward measure, denoted by $f$, is a $\sin^{n-1}$-concave function which attains the maximum at point $o$, thanks to central symmetry of $K_1$. Then we have :
\begin{eqnarray*}
\frac{vol_n((\mathbb{S}^{n-1}_{o}+\varepsilon)\cap K_1)}{vol_n(K_1)} &=& \frac{\int_{[-\varepsilon,\varepsilon]\subset \sigma}f(t)dt}{\int_{\sigma}f(t)dt}\\
                                                                     &\geq& \frac{\int_{0}^{\varepsilon}\cos^{n-1}(t)dt}{\int_{0}^{\pi/2}\cos^{n-1}(t)dt} \\
                                                                     &=& \frac{vol_n(\mathbb{S}^{n-1}+\varepsilon)}{vol_n(\mathbb{S}^n)},
\end{eqnarray*}
where the two last equations are obtained by applying Lemma \ref{grmem}.

This completes the proof of Proposition \ref{import}.

\begin{flushright}
$\Box$
\end{flushright}

\emph{Remark}: Generalising Lemma \ref{brun} to higher dimensions and using Lemma \ref{grmem}, we can thus easily obtain the following:
\begin{prop} \label{importt}
For every $\varepsilon>0$ and every $\mathbb{S}^{n-k}_{o}$ we have :
\begin{eqnarray*}
\frac{vol_n((\mathbb{S}^{n-k}_{o}+\varepsilon)\cap K_1)}{vol_n(K_1)}\geq \frac{vol_n(\mathbb{S}^{n-k}_{o}+\varepsilon)}{vol_n(\mathbb{S}^n)}.
\end{eqnarray*}
\end{prop}

This proves that every symmetric convex set and every tube of the form $\mathbb{S}^k+\varepsilon$ satisfy the correlation inequality of Conjecture \ref{princ}.

Lemmas \ref{bg} and \ref{high} combined with propositions \ref{equal} and \ref{importt} complete the proof of Theorem \ref{dmain}.

\begin{flushright}
$\Box$
\end{flushright}

\section{From the Euclidean Space to the Sphere and the Cauchy Correlation Conjecture}

In this section, we shall prove that Conjecture \ref{pri} is equivalent to Conjecture \ref{princ}. For this we shall recall an important map between the Euclidean space to the open hemi-sphere:

\subsection{The Projective (or Gnomonic Projection) of the Euclidean Space onto the Sphere}

Let $\mathbb{R}^n=\mathbb{R}^n\times{0}\subset \mathbb{R}^{n+1}$. Let $\mathbb{S}^n_{-1}$ be the unit sphere centered at $c=(0,\cdots,-1)$. 
\begin{de}[Gnomonic Projection] \label{gno}
The Gnomonic map is the map $q: \mathbb{R}^n\rightarrow \mathbb{S}^n_{-1}$ defined such that for every $x\in \mathbb{R}^n$
\begin{eqnarray*}
q(x)=[x,c]\cap \mathbb{S}^n_{-1},
\end{eqnarray*}
where $[x,c]$ is the \emph{line} segment (in $\mathbb{R}^{n+1}$) joining the center of the sphere to the point $x$.
\end{de}

$p$ is an isomorphism between the Euclidean space and a (fixed) ball of radius $\pi/2$ of the sphere. By the previous definition, one can clearly see that $p$ maps every straight line to a geodesic segment. This particularity of the gnomonic projection is essential for us since $p$ maps convex sets of Euclidean space to convex subsets of a hemi-sphere of the round sphere. Note that  the map $p$ is neither the exponential map (for which only the lines passing through the origin are mapped to geodesics) nor the stereographical map (for which the image of a point of the sphere is the intersection of a segment passing through this point and a fixed point of the sphere with the Euclidean space). This specific map was used in \cite{gromil} and \cite{grwst} to transport measures supported by convex subsets of the Euclidean space to measures supported by convex sets of the sphere and is used to prove some isoperimetric type inequalities on the sphere and the Gaussian space.
\begin{lem} \label{2}
The push-forward of the Cauchy measure $\nu_n$ with density 
\begin{eqnarray*}
p(x)=C.(1+\vert x\vert^2)^{-\frac{n+1}{2}}
\end{eqnarray*}
 under the gnomonic projection is the normalised canonical Riemannian measure $\frac{1}{vol(B(.,\pi/2))}dv_{\mathbb{S}^n_{+}}$.
\end{lem}

\emph{Proof of the Lemma \ref{2}} 

It is sufficient to calculate the Jacobian of the map of definition \ref{gno}.

\begin{flushright}
$\Box$
\end{flushright}

\subsection{Equivalence Between Conjecture \ref{princ} and Conjecture \ref{pri}}

\begin{prop} \label{eqi}
Conjecture \ref{pri} is equivalent to Conjecture \ref{princ}.
\end{prop}

\emph{Proof of Proposition \ref{eqi}}



Applying the isomorphism $q$ of definition \ref{gno} and Lemma \ref{2}, we can translate/transport the setting of Conjecture \ref{pri} on a hemi-sphere of $\mathbb{S}^n$ :

Let $K$ and $M$ be the two convex bodies of the assumption of Conjecture \ref{pri}. Let $K_1=q(K)$ and $K_2=q(M)$. Applying Lemma \ref{2} one has 
\begin{eqnarray*}
\nu_n(K)&=&q^{-1}_{*}(\mu_n)(K_1)=vol_n(K_1)/vol_n(B(o,\pi/2)) \\ 
\nu_n(M)&=&q^{-1}_{*}(\mu_n)(K_2)=vol_n(K_2)/vol_n(B(o,\pi/2)) \\
\nu_n(K\cap M)&=&q^{-1}_{*}(\mu_n)(K_1\cap K_2)=vol_n(K_1\cap K_2)/vol_n(B(o,\pi/2)).
\end{eqnarray*}
Note that $q(K\cap M)=q(K)\cap q(M)$ since of course $q$ is an isomorphism. If Conjecture \ref{princ} holds, we get
\begin{eqnarray*}
\nu_n(K\cap M)&=&\frac{vol_n(K_1\cap K_2).vol_n(B(o,\pi/2))}{(vol_n(B(o,\pi/2)))^2} \\
              &\geq& \frac{(vol_n(K_1)}{vol_n(B(o,\pi/2))}.\frac{vol_n(K_2)}{vol_n(B(o,\pi/2))} \\
              &=&\nu_n(K).\nu_n(M).
\end{eqnarray*}
This proves a correlation theorem for Cauchy Measures.

If one supposes Conjecture \ref{pri} holds (again by applying the above argument using the Gnomonic projection) we can obtain the spherical correlation and hence Conjecture \ref{princ}.

This ends the proof of Proposition \ref{eqi}.

\begin{flushright}
$\Box$
\end{flushright}


\section{Remarks and Questions}

\begin{itemize}

\item Is there an $O(n)$-invariant map from $\mathbb{R}^n$ to $\mathbb{S}^n_{+}$ which transports the Gaussian measure to the normalised Riemannian measure of the hemi-sphere (such that given two (bounded) convex subsets $K_1$ and $K_2$ of $\mathbb{R}^n$, their image remain convex)?

If such a map exists, then the Gaussian Correlation Conjecture can be proved directly from Conjecture \ref{princ} using this map.

\item Following the proof of Conjecture \ref{princ} for the special cases discussed in Section $4$, it seems that the spherical correlation Conjecture should hold for a wider class of sets. My first supposition is that Conjecture \ref{princ} should hold for a pair of symmetric convex set and a symmetric \emph{mean-convex} set (a set with positive mean-curvature of boundary). I did not follow up on this problem, but it could be interesting to characterise every two subsets (not necessarely convex) of the sphere which satisfy the correlation property.

\item One possible way of proving Conjecture \ref{princ} (according to Section $4$) would be, for example, to choose a \emph{good} $k$-dimensional section of one set $K_1$, replace it by a tube $\mathbb{S}^k+\varepsilon$ of the same volume, deform the other convex set $K_2$ to become a generalised double cone of appropriate (co)-dimension and same volume. Then, according to Proposition \ref{equal}, we have the equality case for the spherical correlation and it would remain to prove that the volume of the intersection of this deformation becomes smaller than the original volume of the intersection. I attempted this method by trying to find the \emph{good} section  by generalising the Dvoretzky Theorem for two symmetric convex sets but still wasn't able to prove Conjecture \ref{princ} in its most general form.

\item In the past few years localisation methods were used to prove very interesting geometric inequalities. In \cite{lova} and \cite{kann} the authors prove integral formulae using localisation, and apply their methods to conclude a few isoperimetric type inequalities concerning the convex sets in the Euclidean space. In \cite{guedon} the authors study a functional analysis version of the localisation used again on the Euclidean space. The localisation on more general spaces was studied in \cite{gromil}, \cite{grwst}, \cite{memwst}, \cite{memusphere} and \cite{membru}. It may seem hard to believe that one could prove Conjecture \ref{princ} using localisation methods but an easier version of Conjecture \ref{princ}, where we replace $vol_n(B(o,\pi/2)$ by $vol_n(\mathbb{S}^n)$, should be possible to prove using the following proposition proved in \cite{membru} :
\begin{prop} \label{fonda}
Let $f_1$, $f_2$ be two upper semi-continuous nonegative functions on $\mathbb{S}^n$ and $f_3$, $f_4$ be two lower semi-continuous nonegative functions on $\mathbb{S}^n$. Let $-\infty\leq s\leq 1/2$ and $\alpha$,$\beta>0$. Suppose that $f_1^{\alpha}f_2^{\beta}\leq f_3^{\alpha}f_4^{\beta}$ and for every $a,b\in \mathbb{S}^n$, for every $\sin^s$-affine probability measure $\nu$ supported by the geodesic segment $[a,b]$,
\begin{eqnarray*}
(\int f_1 d\nu)^{\alpha}(\int f_2 d\nu)^{\beta}\leq (\int f_3 d\nu)^{\alpha}(\int f_4 d\nu)^{\beta}.
\end{eqnarray*}
Then for every $\sin^s$-concave probability measure $\mu$ on $\mathbb{S}^n$,
\begin{eqnarray*}
(\int f_1 d\mu)^{\alpha}(\int f_2 d\mu)^{\beta}\leq (\int f_3 d\mu)^{\alpha}(\int f_4 d\mu)^{\beta}.
\end{eqnarray*}
\end{prop}
To prove Conjecture \ref{princ} applying Proposition \ref{fonda}, one can fix $K_1$ and translate (without rotating) $K_2$ to the point $-o$. Distribute the mass of $K_1\cap K_2$ around the waist of the sphere (i.e. around $\partial B(o,\pi/2)=\mathbb{S}^{n-1}$). I believe this is a suitable geometric configuration for one to be able to apply Proposition \ref{fonda} and prove a weak version of Conjecture \ref{princ}. 

\item The Bishop-Gromov Inequality  asserts that for any convex set $X\subset \mathbb{S}^n$ and for every $x\in X$, the function
\begin{eqnarray*}
\frac{vol_n(X\cap B(x,r))}{vol_n(B(x,r))}
\end{eqnarray*}
is a non-increasing function of $r$ (See Lemma \ref{bg}). If the following conjecture (which is a generalisation of the Bishop-Gromov Inequality) holds, then the proof of Conjecture \ref{princ} becomes straightforward.
\begin{conj}  \label{bgr}
Let $X$ be a convex subset of the sphere having $o\in\mathbb{S}^n$ in its interior. Let $F_1\subset F_2$ be two symmetric convex subsets of the sphere, both containing $o$ in their interiors. Then
\begin{eqnarray*}
\frac{vol_n(X\cap F_1)}{vol_n(F_1)}\geq \frac{vol_n(X\cap F_2)}{vol_n(F_2)}
\end{eqnarray*}
\end{conj}
I attempted to prove this conjecture by putting a density on the sphere and using the Bishop-Gromov Inequality for Riemannian manifolds with density having positive Ricci curvature in the sense of Lott-Villani (see \cite{lotvil}) but I was not successful. A possible way to prove Conjecture \ref{bgr} could be the use of the theory of mean-curvature flow (see \cite{mf}).

 
\item The theory of Ricci curvature for general metric-measure spaces is well developed. For an example, one can see \cite{lotvil} where the authors prove Bishop-Gromov type inequality for metric-measure spaces having $Ricci>0$ in the sense of displacement convexity. By using their definition of  Ricci curvature and combining it with their Bishop-Gromov type inequality for the Gaussian space, can one directly (without passing through the sphere) prove Conjecture \ref{pri}?

\item What type of correlation theorem (like the one we proposed in Conjecture \ref{princ}) can one prove for Riemannian manifolds with a lower bound on the Ricci curvature?
\end{itemize}

\bibliographystyle{plain}
\bibliography{gcor}

\begin{thebibliography}{10}

\bibitem{almg}
F.~{Almgren}.
\newblock Optimal isoperimetric inequalities.
\newblock {\em Indiana University Mathematics Journal}, 35(3):451--547, 1986.

\bibitem{dasgu}
S.~{Das Gupta}, M.L. {Eaton}, I.~{Olkin}, M.~{Perlman}, L.J. {Savage}, and
  M.~{Sobel}.
\newblock Inequalities on the probability content of convex regions for
  elliptically contoured distributions.
\newblock {\em proc. Sixth Berkeley Symp. Math. Statist. Prob.}, 3:241--264,
  1972.

\bibitem{dunnet}
C.W. {Dunnet} and M.~{Sobel}.
\newblock Approximations to the probability integral and certain percentage
  points to a multivariate analogue of {S}tudent's t-distribution.
\newblock {\em Biometrika}, 42:258--260--405, 1955.

\bibitem{figa}
A.~{Figalli}, F.~{Maggi}, and A.~{Pratelli}.
\newblock A geometric approach to correlation inequalities in the plane.
\newblock {\em Ann.Ins.H.Poinacr\'e}, to appear.

\bibitem{guedon}
M.~{Fradelizi} and O.~{Gu\'edon}.
\newblock The extreme points of subsets of s-concave probabilities and a
  geometric localization theorem.
\newblock {\em Discrete Comput. Geom.}, 31:327--335, 2004.

\bibitem{membru}
M.~{Fradelizi}, O.~{Gu\'edon}, and Y.~{Memarian}.
\newblock A {B}runn-{M}inkowski type inequality on the sphere.
\newblock {\em preprint}, 2011.

\bibitem{grwst}
M.~{Gromov}.
\newblock Isoperimetry of waists and concentration of maps.
\newblock {\em GAFA}, 13:178--215, 2003.

\bibitem{groexp}
M.~{Gromov}.
\newblock Singularities, expanders and topology of maps. part2: From
  combinatorics to topology via algebraic isoperimetry.
\newblock {\em Geom. Funct. Anal}, 20:416--526, 2010.

\bibitem{gromil}
M.~{Gromov} and V.D. {Milman}.
\newblock Generalisation of the spherical isoperimetric inequality to uniformly
  convex {B}anach spaces.
\newblock {\em Compositio Math.}, 62:3:263--282, 1987.

\bibitem{kann}
R.~{Kannan}, L.~{Lov\'asz}, and M.~{Simonovits}.
\newblock Isoperimetric problems for convex bodies and a localization lemma.
\newblock {\em Discrete Comput. Geom.}, 13(3--4):541--559, 1995.

\bibitem{lewi1}
T.~{Lewis} and G.~{Pritchard}.
\newblock Correlation measures.
\newblock {\em Elec. Jour. Prob.}, 4:77--85, 1999.

\bibitem{lewi2}
T.~{Lewis} and G.~{Pritchard}.
\newblock Tail properties of correlation measures.
\newblock {\em Jour. Theoret. Probab.}, 16:771--788, 2003.

\bibitem{lotvil}
J.~{Lott} and C.~{Villani}.
\newblock Ricci curvature for metric-measure spaces via optimal transport.
\newblock {\em Annals of {M}athematics}, 169:903--991, 2009.

\bibitem{lova}
L.~{Lov\'asz} and M.~{Simonovits}.
\newblock Random walks in a convex body and an improved volume algorithm.
\newblock {\em Random {S}tructures {A}lgorithms}, 4(4):359--412, 1993.

\bibitem{mempos}
Y.~{Memarian}.
\newblock A note on the geometry of positively curved {R}iemannian manifolds.
\newblock {\em In preparation}.

\bibitem{memwst}
Y.~{Memarian}.
\newblock On {G}romov's waist of the sphere theorem.
\newblock {\em Journal of {T}opology and {A}nalysis}, 3:7--36, 2011.

\bibitem{memusphere}
Y.~{Memarian}.
\newblock A lower bound on the waist of unit spheres of uniformly normed
  spaces.
\newblock {\em Compositio Math.}, 148(4):1238--1264, 2012.

\bibitem{mf}
X.P. {Zhu}.
\newblock {\em Lectures on Mean Curvature Flows}.
\newblock AMS/IP, 2002.

\end{thebibliography}

\end{document}